\documentclass{article}
\usepackage{amssymb}
\usepackage{amsmath,color}
\usepackage[hidelinks]{hyperref}

\usepackage[square,numbers]{natbib}         
\newtheorem{lem}{Lemma}[section]
\newtheorem{cor}[lem]{Corollary}
\newtheorem{teo}[lem]{Theorem}
\newtheorem{os}[lem]{Remark}

\newtheorem{prop}[lem]{Proposition}

\newcommand{\qed}{\thinspace\null\nobreak\hfill\hbox{\vbox{\kern-.2pt\hrule
 height.2pt depth.2pt\kern-.2pt\kern-.2pt \hbox to2.5mm{\kern-.2pt\vrule
 width.4pt \kern-.2pt\raise2.5mm\vbox to.2pt{}\lower0pt\vtop
 to.2pt{}\hfil\kern-.2pt \vrule
 width.4pt \kern-.2pt}\kern-.2pt\kern-.2pt\hrule height.2pt depth.2pt
 \kern-.2pt}}\par\medbreak}

\newcommand{\R}{\mathbb{R}}

\newcommand{\C}{\mathbb{C}}

\newcommand{\N}{\mathbb{N}}

\newcommand{\eps}{\varepsilon}

\newcommand{\ds}{\displaystyle}

\date{}
\begin{document}

\title{Anisotropic Sobolev spaces with weights}
\author{G. Metafune \thanks{Dipartimento di Matematica e Fisica ``Ennio De Giorgi'', Universit\`a del Salento, C.P.193, 73100, Lecce, Italy.
-mail:  giorgio.metafune@unisalento.it}\qquad L. Negro \thanks{Dipartimento di Matematica e Fisica  ``Ennio De
Giorgi'', Universit\`a del Salento, C.P.193, 73100, Lecce, Italy. email: luigi.negro@unisalento.it} \qquad C. Spina \thanks{Dipartimento di Matematica e Fisica``Ennio De Giorgi'', Universit\`a del Salento, C.P.193, 73100, Lecce, Italy.
e-mail:  chiara.spina@unisalento.it}}

\maketitle
\begin{abstract}
\noindent 
We study Sobolev spaces with weights in the half-space $\R^{N+1}_+=\{(x,y): x \in \R^N, y>0\}$, adapted to the   singular elliptic   operators 
 \begin{equation*}
\mathcal L =y^{\alpha_1}\Delta_{x} +y^{\alpha_2}\left(D_{yy}+\frac{c}{y}D_y  -\frac{b}{y^2}\right).
\end{equation*}

\par

\end{abstract}

\section{Introduction}
Elliptic and parabolic problems associated to the  degenerate   operators 
\begin{equation*} \label{defL}
\mathcal L =y^{\alpha_1}\Delta_{x} +y^{\alpha_2}\left(D_{yy}+\frac{c}{y}D_y  -\frac{b}{y^2}\right) \quad {\rm and}\quad D_t- \mathcal L
\end{equation*}
in the half-space $\R^{N+1}_+=\{(x,y): x \in \R^N, y>0\}$ or  in $(0, \infty) \times \R^{N+1}_+$ lead quite naturally to the introduction of  weighted Sobolev spaces which are  anisotropic if $\alpha_1 \neq \alpha_2$. The aim of this paper is to provide the functional analytic properties of these Sobolev spaces needed in \cite{MNS-CompleteDegenerate} and in \cite{MNS-PerturbedBessel} in the $1$-d case, where we prove existence, uniqueness and regularity of elliptic and parabolic problems governed by the operators above. We  also refer to   \cite{met-calv-negro-spina, MNS-Sharp, MNS-Grad,  MNS-Grushin, MNS-Max-Reg,  Negro-Spina-Asympt} for the analogous results concerning the $N$-d version of $D_{yy}+\frac{c}{y}D_y  -\frac{b}{y^2}$.

For $m \in \R$ we consider the measure $y^m dx dy $ in $\R^{N+1}_+$ and  write $L^p_m$ for  $L^p(\R_+^{N+1}; y^m dx dy)$. 
Given  $p>1$, $\alpha_1 \in \R$, $\alpha_2<2$, we define the Sobolev space
\begin{align*}
W^{2,p}(\alpha_1,\alpha_2,m)&=\left\{u\in W^{2,p}_{loc}(\R^{N+1}_+):\ u,\  y^{\alpha_1} D_{x_ix_j}u,\ y^\frac{\alpha_1}{2} D_{x_i}u,  y^{\alpha_2}D_{yy}u,\ y^{\frac{\alpha_2}{2}}D_{y}u\in L^p_m\right\}
\end{align*}
which is a Banach space equipped with the norm
\begin{align*}
\|u\|_{W^{2,p}(\alpha_1,\alpha_2,m)}=&\|u\|_{L^p_m}+\sum_{i,j=1}^n\|y^{\alpha_1} D_{x_ix_j}u\|_{L^p_m}+\sum_{i=1}^n\|y^{\frac{\alpha_1}2} D_{x_i}u\|_{L^p_m}\\
&+\|y^{\alpha_2}D_{yy}u\|_{L^p_m}+\|y^{\frac{\alpha_2}{2}}D_{y}u\|_{L^p_m}.
\end{align*}
Next we add a Neumann boundary condition for $y=0$  in the form $y^{\alpha_2-1}D_yu\in L^p_m$ and set
\begin{align*}
W^{2,p}_{\mathcal N}(\alpha_1,\alpha_2,m)=\{u \in W^{2,p}(\alpha_1,\alpha_2,m):\  y^{\alpha_2-1}D_yu\ \in L^p_m\}
\end{align*}
with the norm
$$
\|u\|_{W^{2,p}_{\mathcal N}(\alpha_1,\alpha_2,m)}=\|u\|_{W^{2,p}(\alpha_1,\alpha_2,m)}+\|y^{\alpha_2-1}D_yu\|_{ L^p_m}.
$$

We  consider also an integral version of the Dirichlet boundary condition, namely  a weighted summability requirement for $y^{-2}u$ and introduce 
$$
	W^{2,p}_{\mathcal R}(\alpha_1, \alpha_2, m)=\{u \in  W^{2,p}(\alpha_1, \alpha_2, m): y^{\alpha_2-2}u \in L^p_m\}
$$
with the norm $$\|u\|_{W^{2,p}_{\mathcal R}(\alpha_1, \alpha_2, m)}=\|u\|_{W^{2,p}(\alpha_1, \alpha_2, m)}+\|y^{\alpha_2-2}u\|_{L^p_m}.$$

Note that $\alpha_1, \alpha_2$ are not assumed to be positive. The restriction $\alpha_2<2$ is not really essential since one can deduce from it the case $\alpha_2>2$, using the change of variables described in the next section. However, we keep it both to simplify the exposition and because $\mathcal L$ is mainly considered for $\alpha_2<2$. 

No requirement is made for the mixed derivatives $D_{x_iy}u$ to simplify some arguments. However, the weighted integrability of the mixed derivatives is automatic under the condition of Proposition \ref{Sec sob derivata mista}.

Sobolev spaces with weights are well-known in the literature, see e.g. \cite{grisvard}, \cite[Chapter 6]{necas}, \cite{Geymonat-Grisvard} and \cite{morel} for the non-anisotropic case. Variants of $W^{2,p}_{\mathcal R}(\alpha_1, \alpha_2, m)$, usually defined as the closure of compactly supported functions in $W^{2,p}(\alpha_1, \alpha_2, m)$, can be found in the above papers . However, we have not been able to find anything about  $W^{2,p}_{\mathcal N}(\alpha_1, \alpha_2, m)$.

Let us briefly describe the content of the paper.  In Section 2 we show that, by a change of variables, the spaces 
$W^{2,p}_{\mathcal{N}}(\alpha_1,\alpha_2,m)$ and $W^{2,p}_{\mathcal{N}}(\tilde \alpha_1,\tilde\alpha_2,\tilde m)$, $ \tilde\alpha_1=\frac{\alpha_1}{\beta+1},\quad \tilde\alpha_2=\frac{\alpha_2+2\beta}{\beta+1}, \quad \tilde m=\frac{m-\beta}{\beta+1}$ are isomorphic. This observation simplifies many proofs but requires the full scale of $L^p_m$ spaces, according to the general strategy of \cite{} to study the operator $\mathcal L$. Hardy inequalities and traces for $y=0$ are studied in Section 3. The main properties of the spaces $W^{2,p}_{\mathcal{N}}(\alpha_1,\alpha_2,m)$ are proved in Section 4 together with a density result for smooth functions having zero $y$-derivative in a strip around $y=0$, which  is crucial in the study of the operator $\mathcal L$. The space $W^{2,p}_{\mathcal{R}}(\alpha_1,\alpha_2,m)$ is studied in Section 5.

\section{A useful change of variables}\label{Section Degenerate}

For  $k,\beta \in\R$, $\beta\neq -1$ let
\begin{align}\label{Gen Kelvin def}
T_{k,\beta\,}u(x,y)&:=|\beta+1|^{\frac 1 p}y^ku(x,y^{\beta+1}),\quad (x,y)\in\R^{N+1}_+.
\end{align}
Observe that
$$ T_{k,\beta\,}^{-1}=T_{-\frac{k}{\beta+1},-\frac{\beta}{\beta+1}\,}.$$

\begin{lem}\label{Isometry action der} The following properties hold for $1 \leq p \leq \infty$.
	\begin{itemize}
		\item[(i)] $T_{k,\beta\,}$ maps isometrically  $L^p_{\tilde m}$ onto $L^p_m$  where 
		$$ \tilde m=\frac{m+kp-\beta}{\beta+1}.$$
		\item[(ii)] For every  $u\in W^{2,1}_{loc}\left(\R^{N+1}_+\right)$ one has
\begin{itemize}
	\item[1.] $y^\alpha T_{k,\beta\,}u=T_{k,\beta\,}(y^{\frac{\alpha}{\beta+1}}u)$, for any $\alpha\in\R$;\medskip
	\item [2.] $D_{x_ix_j}(T_{k,\beta\,}u)=T_{k,\beta} \left(D_{x_ix_j} u\right)$, \quad $D_{x_i}(T_{k,\beta\,}u)=T_{k,\beta}\left(D_{x_i} u\right)$;\medskip
	\item[3.]  $D_y T_{k,\beta\,}u=T_{k,\beta\,}\left(ky^{-\frac 1 {\beta+1}}u+(\beta+1)y^{\frac{\beta}{\beta+1}}D_yu\right)$,
	\\[1ex] $D_{yy} (T_{k,\beta\,} u)=T_{k,\beta\,}\Big((\beta+1)^2y^{\frac{2\beta}{\beta+1}}D_{yy}u+(\beta+1)(2k+\beta)y^{\frac{\beta-1}{\beta+1}}D_y u+k(k-1)y^{-\frac{2}{\beta+1}}u\Big)$.\medskip
	\item[4.] $D_{xy} T_{k,\beta\,}u=T_{k,\beta\,}\left(ky^{-\frac 1 {\beta+1}}D_xu+(\beta+1)y^{\frac{\beta}{\beta+1}}D_{xy}u\right)$
\end{itemize}
	\end{itemize}
\end{lem}{\sc{Proof.}} The proof of (i) follows after observing the Jacobian of $(x,y)\mapsto (x,y^{\beta+1})$ is $|1+\beta|y^{\beta}$. To prove (ii) one first observes that any $x$-derivatives commutes with $T_{k,\beta}$. Then we compute  
\begin{align*}
D_y T_{k,\beta\,}u(x,y)=&|\beta+1|^{\frac 1 p}y^{k}\left(k\frac {u(x,y^{\beta+1})} y+(\beta+1)y^\beta D_y u(x,y^{\beta+1})\right)\\[1ex]
=&T_{k,\beta\,}\left(ky^{-\frac 1 {\beta+1}}u+(\beta+1)y^{\frac{\beta}{\beta+1}}D_yu\right)
\end{align*}
and similarly
\begin{align*}
D_{yy} T_{k,\beta\,} u(x,y)=&T_{k,\beta\,}\Big((\beta+1)^2y^{\frac{2\beta}{\beta+1}}D_{yy}u+(\beta+1)(2k+\beta)y^{\frac{\beta-1}{\beta+1}}D_y u+k(k-1)y^{-\frac{2}{\beta+1}}u\Big).
\end{align*}
\qed

Let us specialize the above lemma to 
\begin{align*}
T_{0,\beta}&:L^p_{\tilde m}\to L^p_m,\qquad  \tilde m=\frac{m-\beta}{\beta+1}
\end{align*}
to transform Sobolev spaces with different exponents.

\begin{prop}
\label{Sobolev eq}
Let $p>1$, $m, \alpha_1,\alpha_2\in \R$ with $\alpha_2< 2$.
 Then one has 
\begin{align*}
W^{2,p}_{\mathcal{N}}(\alpha_1,\alpha_2,m)=T_{0,\beta}\left(W^{2,p}_{\mathcal{N}}(\tilde \alpha_1,\tilde\alpha_2,\tilde m)\right),\qquad \tilde\alpha_1=\frac{\alpha_1}{\beta+1},\quad \tilde\alpha_2=\frac{\alpha_2+2\beta}{\beta+1}.
\end{align*}
In particular, by choosing $\beta=-\frac{\alpha_2}2$ one has 
\begin{align*}
W^{2,p}_{\mathcal{N}}(\alpha_1,\alpha_2,m)=T_{0,-\frac {\alpha_2} 2}\left(W^{2,p}_{\mathcal{N}}(\tilde \alpha,0,\tilde m)\right),\qquad \tilde\alpha=\frac{2\alpha_1}{2-\alpha_2},\quad \tilde m=\frac{m+\frac{\alpha_2} 2}{1-\frac{\alpha_2} 2}.
\end{align*}
\end{prop}
{\sc{Proof.}} Given  $\tilde u\in W^{2,p}_{\mathcal{N}}(\tilde \alpha_1,\tilde\alpha_2,\tilde m)$  let us set $ u(x,y)=(T_{0,\beta}\tilde u)(x,y)=|\beta+1|^{1/p}\tilde u(x,y^{\beta+1})$. Everything follows from the equalities of Lemma \ref{Isometry action der},
\begin{itemize}
\item [(i)] $y^{\alpha_1}D_{x_ix_j}u=T_{0,\beta} \left(y^{\tilde \alpha_1}D_{x_ix_j}\tilde u\right)$, \quad $y^{\frac{\alpha_1}{2}}D_{x_i}u=T_{0,\beta}\left(y^{\frac{\tilde \alpha_1}{2}}D_{x_i}\tilde u\right)$;\smallskip
\item [(ii)] $y^{\frac{\alpha_2}2}D_{y}u=(1+\beta)T_{0,\beta}\left(y^{\tilde\alpha_2}D_{y}\tilde u\right)$,\quad $y^{\alpha_2-1}D_{y}u=(1+\beta)T_{0,\beta}\left(y^{\tilde \alpha_2-1}D_{y}\tilde u\right)$;\smallskip
\item [(iii)] $y^{\alpha_2}D_{yy}u=(1+\beta)T_{0,\beta}\left[(1+\beta) y^{\tilde \alpha_2}D_{yy}\tilde u+\beta y^{\tilde \alpha_2-1}D_{y}\tilde u\right]$.
\end{itemize}
\qed

\begin{os}
	Note that in the above  proposition is essential to deal with  $W^{2,p}_{\mathcal{N}}(\alpha_1,\alpha_2,m)$. Indeed in  general the isometry $T_{0,\beta}$ does  not  transform  $W^{2,p}(\tilde \alpha_1,\tilde\alpha_2,\tilde m)$ into $W^{2,p}(\alpha_1,\alpha_2,m)$, because of identity (iii) above.
\end{os}

\section{Hardy inequalities and traces}
In this section we prove some weighted Hardy inequalities and investigate trace properties of function $u$ such that $y^\beta D_yu\in L^p_m$.

\medskip
The following  result is standard but we give a proof to settle "almost everywhere" issues.

\begin{lem} \label{L1}
	 Let $u\in L^1_{loc}(\R^{N+1}_+)$ be such that $D_yu\in L^1(\R^{N+1}_+)$. Then  there exists $v$ such that  $v=u$ almost everywhere and  $v(\cdot,y)\in L^1_{loc}(\R^N)$ for every $y\geq 0$ and  
	$$v(x,y_2)-v(x,y_1)=\int_{y_1}^{y_2}D_yu(x,s)\ ds$$ for every $0 \leq y_1< y_2 \leq \infty$ and  almost every $x\in Q$.
\end{lem}
{\sc Proof.} For a.e.  $x \in \R^N$ the function $u(x, \cdot)$ is absolutely continuous and then
$$u(x,y_2)-u(x,y_1)=\int_{y_1}^{y_2}D_yu(x,s)\ ds$$
for a.e. $y_1, y_2$. It is therefore sufficient to define $v(x,y)=\int_c^y D_yu(x,s)\, ds+u(x,c)$,  if $c$ is chosen in such a way that $u(\cdot,c) \in L^1_{loc}(\R^N)$.
\qed

Properties of  functions  $u\in L^p_m$ such that $D_y u \in L^p_m$ have been proved in \cite[Appendix B]{MNS-Caffarelli}. Here we exploit the more general property $y^\beta D_y u \in L^p_m$.
%
%

\begin{prop}
	\label{Hardy in core}
	Let  $C:=\left|\frac{m+1}{p}-(1-\beta)\right|^{-1}$. The following properties hold for  $ u\in L^1_{loc}(\R^{N+1}_+)$ such that $y^{\beta}D_{y}u\in L^p_m$. 
\begin{itemize}
	\item[(i)] If  $\frac{m+1}p<1-\beta$  then $D_yu\in L^1\left(Q\times  [0,1]\right)$ for any cube $Q$ of $\R^N$; in particular  $u$ has a trace $u(\cdot,y)\in L^1_{loc}(\R^N)$ for every $0\leq y\leq 1$. Moreover setting $u_{0}(x)=\lim_{y\to0}u(x,y)$ one has 
	\begin{align*}
		\|y^{\beta-1}(u-u_0)\|_{L^p_m}\leq C \|y^{\beta}D_{y}u\|_{L^p_m}.
	\end{align*}
If  moreover $u\in L^p_m$ then $u(\cdot,y)\in L^p(\R^N)$ for every $0\leq y\leq 1$.	
	\item[(ii)] If $\frac{m+1}p>1-\beta$  then $D_yu\in L^1\left(Q\times  [1,\infty[\right)$ for any cube $Q$ of $\R^N$; in particular $u$  has a finite trace $u_{\infty}(x)=\lim_{y\to\infty}u(x,y)\in L^1_{loc}\left(\R^N\right)$ and 
	\begin{align}
		\|y^{\beta-1}(u-u_{\infty)}\|_{L^p_m}\leq C \|y^{\beta}D_{y}u\|_{L^p_m}.
	\end{align}
If moreover $u\in L^p_m$ then $u_\infty\in L^p(\R^N)$ and $u_\infty=0$ if $m \geq -1$.
	
\end{itemize}	
\end{prop}
{\sc{Proof.}} To prove (i) let  $f(x,y):=y^{\beta}D_{y}u(x,y)$. If  $Q$ is a cube of $\R^N$ then since $\frac{m+1}p>1-\beta$ one has 
\begin{align*}
	\int_{Q\times [0,1]}|D_yu|dxdy&=\int_{Q\times [0,1]}|D_yu|y^{\beta}y^{-\beta-m}y^mdxdy\\[1ex]
	&\leq \|y^{\beta}D_yu\|_{L^p_m}\left(\int_{0}^1 y^{-(\beta+m)p'+m}\right)^{\frac 1{p'}}|Q|^{\frac 1{p'}}
	=C(Q,b,p)\|y^{\beta}D_yu\|_{L^p_m}.
\end{align*}
In particular by Lemma \ref{L1},  $u$ has a finite trace $u(\cdot,y)\in L^1_{loc}\left(\R^N\right)$
 for every $0\leq y\leq 1$. Setting
 $u_0(x)=u(x,0)=\lim_{y \to 0}u(x,y ) $ we can write
\begin{align*}
	y^{\beta-1}\left(u(x,y)-u_{0}(x)\right)=y^{\beta-1}\int_{0}^y f(x,s)s^{-\beta}\,ds:=(H_1f)(y).
\end{align*}
By \cite[Lemma 10.3, (i)]{MNS-Caffarelli}, the operator $H_1$ is bounded on $L^p_m(\R_+)$ when $\frac{m+1}p<1-\beta$, hence
\begin{align*}
	\|y^{\beta-1}\left (u(x,\cdot)-u_0 (x) \right)\|_{L^p_m(\R_+)}\leq C \|y^{\beta}D_{y}u(x,\cdot)\|_{L^p_m(\R_+)}.
\end{align*} 
Claim (i) then  follows by raising to the power $p$ and integrating with respect to $x$. To prove that $u(\cdot,y)\in L^p(\R^N)$ we proceed analogously: since $u\in L^p_m$ then $u(\cdot,y)\in L^p(\R^N)$ for a.e. $y\in [0,1]$. Without any  loss of generality we suppose $u(\cdot,1)\in L^p(\R^N)$  and we  write for any $y_0\in [0,1]$
\begin{equation*}
	u(x,y_0)=u(x,1)-\int_{y_0}^ 1D_y u(x,s)\ ds=u(x,1)-\int_{y_0}^1 s^\beta D_y u(x,s)s^{\frac m p}s^{-\beta-\frac m p}\ ds.
\end{equation*} 
Then using H\^older inequality
\begin{align*}
	|u(x,y_0)|&\leq |u(x,1)|+\left(\int_{y_0}^1 \left|s^\beta D_y u(x,s)\right|^p s^{m}\ ds\right)^{\frac 1 p}\left(\int_{y_0}^1 s^{(-\beta-\frac mp)p'}\ ds\right)^{\frac 1 {p'}}\\[1ex]
	&\leq |u(x,1)|+C\left(\int_{0}^1 \left|s^\beta D_y u(x,s)\right|^p s^{m}\ ds\right)^{\frac 1 p}.
\end{align*}
Raising to the power $p$ and integrating with respect to $x$ we obtain
\begin{align*}
	\|u(\cdot,y_0)\|_{L^{p}(\R^N)}\leq C\left(\|u(\cdot,1)\|_{L^p(\R^N)}+\left\|y^\beta D_yu\right\|_{L^p_m}\right).
\end{align*}

The proof of (ii) is similar writing
\begin{align*}
	y^{\beta-1}\left(u(x,y)-u_{\infty}(x)\right)=-y^{\beta-1}\int_{y}^\infty f(x,s)s^{-\beta}\,ds:=-(H_2f)(y)
\end{align*}
and applying \cite[Lemma 10.3, (ii)]{MNS-Caffarelli}. If $u \in L^p_m$ and $m \geq -1$, then $|u(x, \cdot)|^p$ is not summable with respect to $y^m\, dy$ for every $x$ where $u_\infty (x) \neq 0$, hence $u_\infty=0$ a.e.
\qed
 
%

\medskip

\medskip

	In the next lemma we show that $u$ has  has a logarithmic singularity for $y\to 0, \infty$,  when $\frac{m+1}{p}=1-\beta $.
	
	\begin{lem} \label{int-uMaggiore}
		If  $\frac{m+1}{p}=1-\beta $ and  $u, y^{\beta}D_{y}u\in L^p_m$, then 
			\begin{equation} \label{behaviour}
				\left(\int_{\R^N}|u(x,y)|^p\, dx\right)^{\frac 1 p}\leq  \|u(\cdot,1)\|_{L^p(\R^N)}+|\log y|^{\frac{1}{p'}}\|y^{\beta} D_y\|_{L^p_m}.
			\end{equation}
	\end{lem}
	{\sc Proof.} 
	Let $\frac{m+1}{p}=1-\beta$ and set $f=y^\beta D_y\in L^p_m$. Then  for $y\in (0,1)$ one has 
	\begin{align*}
		u(x,y)&=u(x,1)-\int_y^1 D_y u(x,s)\ ds=u(x,1)-\int_y^1s^{-\beta}f(x,s)\ ds\\[1ex]
		&=u(x,1)-\int_y^1s^{-\beta-m} f(x,s)s^m\ ds.
	\end{align*} 
	Therefore, since $(-\beta-m)p'+m=-1$, H\"older inequality yields
	\begin{align*}
		|u(x,y)|&\leq |u(x,1)|+\left(\int_y^1 s^{(-\beta-m)p'}s^m\ ds\right)^\frac{1}{p'}\left(\int_y^1 |f(x,s)|^ps^m\ ds\right)^\frac{1}{p}
		\\[1ex]
		&\leq | u(x,1)|+ |\log y|^\frac{1}{p'}\|f(x,\cdot)\|_{L^p\left((0,1),y^mdy\right)}.
	\end{align*} 
The inequality for $y>1$ is similar.

Since  $u\in L^p_m$ then, as in Proposition \ref{Hardy in core}, we can suppose $u(\cdot,1)\in L^p(\R^N)$ and raising to the power $p$ and integrating with respect to $x$  we conclude the proof. 
		\qed
		
	We also need some elementary interpolative inequalities; the first generalizes \cite[Lemma 4.3]{met-soba-spi-Rellich}  (see also \cite{met-negro-soba-spina}).
		\begin{lem} \label{inter} For $m, \beta \in \R$, $1<p<\infty$ there exist $C>0, \eps_0>0$ such that for every  $u \in W^{2,p}_{loc}((0, \infty))$, $0<\eps <\eps_0$, 
			$$
			\|y^{\beta-1} u'\|_{L^p_m (\R_+)} \leq C \left (\eps \|y^\beta u''\|_{L^p_m(\R_+)} +\frac{1}{\eps} \|y^{\beta-2}u\|_{L^p_m(\R_+)} \right ).
			$$
		\end{lem}
		{\sc Proof. } Changing $\beta$ we may assume that $m=0$. We use the elementary inequality
		\begin{equation} \label{i1}
			\int_a^b |u'(y)|^p\, dy \leq C\left (\eps^p (b-a)^p \int_a^b |u''(y)|^p\, dy+\frac{1}{\eps^p (b-a)^p}\int_a^b |u(y)|^p\, dy\right )
		\end{equation}
		for $\eps \leq \eps_0$, where $\eps_0, C$ are the same as for the unit interval (this follows by scaling). We apply this inequality to each interval $I_n=[2^n, 2^{n+1}[$, $n \in \mathbb Z$
		and multiply by $2^{n(\beta-1)p}$ thus obtaining since $y \approx 2^n$ in $I_n$ 
		$$
		\int_{I_n}y^{(\beta-1)p} |u'(y)|^p\, dy \leq \tilde C\left (\eps^p \int_ {I_n}y^{\beta p}|u''(y)|^p\, dy+\frac{1}{\eps^p}\int_{I_n} y^{(\beta-2)p}|u(y)|^p\, dy\right ).
		$$
		The thesis follows summing over $n$. \qed
		\begin{lem} \label{inter1} For $m, \beta <2$, $1<p<\infty$ there exist $C>0, \eps_0>0$ such that for every  $u \in W^{2,p}_{loc}((1, \infty))$, $0<\eps <\eps_0$, 
			$$
			\|y^{\frac{\beta}{2}} u'\|_{L^p_m((1, \infty))} \leq C \left (\eps \|y^\beta u''\|_{L^p_m((1, \infty))} +\frac{1}{\eps} \|u\|_{L^p_m((1, \infty))} \right ).
			$$
		\end{lem}
		{\sc Proof. }We use \eqref{i1} in $(a_n, a_{n+1})$ where   $a_n=n^{1+\frac{\gamma}{2}}$, so that $a_{n+1}-a_n \approx n^{\frac{\gamma}{2}}$. We multiply both sides by $n^{(1+\frac{\gamma}{2})(m+\frac{\beta p}{2})} \approx y^{m+\frac{\beta p}{2}}$ in $(a_n, a_{n+1})$ and sum over $n$. Choosing $\gamma \geq 0$ in such a way that $\beta=\frac{2\gamma}{2+\gamma}$, the thesis follows.
		\qed

\section{The space  $W^{2,p}_{\mathcal N}(\alpha_1,\alpha_2,m)$} \label{section sobolev}
Let $p>1$, $m, \alpha_1 \in \R$, $\alpha_2<2$. We recall that  
\begin{align*} W^{2,p}_{\mathcal N}(\alpha_1,\alpha_2,m)=\{u \in W^{2,p}(\alpha_1,\alpha_2,m):\  y^{\alpha_2-1}D_yu\ \in L^p_m\}
\end{align*}
with the norm
$$
\|u\|_{W^{2,p}_{\mathcal N}(\alpha_1,\alpha_2,m)}=\|u\|_{W^{2,p}(\alpha_1,\alpha_2,m)}+\|y^{\alpha_2-1}D_yu\|_{ L^p_m}.
$$

We have made the choice  not to  include the mixed derivatives in the definition of  $W^{2,p}_{\mathcal{N}}\left(\alpha_1,\alpha_2,m\right)$ to  simplify some arguments. 
However the following result holds in a range of parameters which is sufficient for the study of the operator $\mathcal L$.

\begin{prop}\label{Sec sob derivata mista} If   $\alpha_2-\alpha_1<2$ and $\alpha_1^{-} <\frac{m+1}p$ then there exists $C>0$ such that 
$$\|y^\frac{\alpha_1+\alpha_2}{2} D_{y}\nabla_x u \|_{ L^p_m} \leq C \|u\|_{W^{2,p}_{\cal N}(\alpha_1, \alpha_2, m)}$$ for every $u \in W^{2,p}_{\mathcal{N}}\left(\alpha_1,\alpha_2,m\right)$.  
\end{prop}
{\sc Proof.}
This follows from \cite[Theorem 7.1]{MNS-CompleteDegenerate}, choosing $c$ sufficiently large therein, so that  $\alpha_1^{-} <\frac{m+1}p<c+1-\alpha_2$.
  \qed

\begin{os}\label{Os Sob 1-d}
	With obvious changes we may consider also the analogous Sobolev spaces on $\R_+$, $W^{2,p}(\alpha_2,m)$ and $W^{2,p}_{\cal N}(\alpha_2, m)$. For example  we have
			$$W^{2,p}_{\mathcal N}(\alpha,m)=\left\{u\in W^{2,p}_{loc}(\R_+):\ u,\    y^{\alpha}D_{yy}u,\ y^{\frac{\alpha}{2}}D_{y}u,\ y^{\alpha-1}D_{y}u\in L^p_m\right\}.$$
		For brevity sake, we consider in what follows,  only the Sobolev spaces on $\R^{N+1}_+$ but all the results of this section will be valid also in $\R_+$ changing the condition $\alpha_1^{-} <\frac{m+1}p$  (which appears in  Sections \ref{denso},  \ref{Sec sob  min domain}) to $0<\frac{m+1}p$.
\end{os}

We  clarify in which sense the condition $y^{\alpha_2-1}D_y u \in L^p_m$ is a Neumann boundary condition.

\begin{prop} \label{neumann} The following assertions hold.
\begin{itemize} 
\item[(i)] If $\frac{m+1}{p} >1-\alpha_2$, then $W^{2,p}_{\mathcal N}(\alpha_1, \alpha_2, m)=W^{2,p}(\alpha_1, \alpha_2, m)$.
\item[(ii)] If $\frac{m+1}{p} <1-\alpha_2$, then $$W^{2,p}_{\mathcal N}(\alpha_1, \alpha_2, m)=\{u \in W^{2,p}(\alpha_1, \alpha_2, m): \lim_{y \to 0}D_yu(x,y)=0\ {\rm for\ a.e.\   x \in \R^N }\}.$$
\end{itemize}
In both cases (i) and (ii), the norm of $W^{2,p}_{\mathcal N}(\alpha_1, \alpha_2, m)$ is equivalent to that of $W^{2,p}(\alpha_1, \alpha_2, m)$.
\end{prop}
{\sc Proof. } If $\frac{m+1}{p} >1-\alpha_2$  and $u \in W^{2,p}(\alpha_1, \alpha_2, m)$, we apply Proposition  \ref{Hardy in core} (ii) to $D_y u$ and obtain that $\lim_{y \to \infty}D_yu(x,y)=g(x)$ exists. At the points where $g(x) \neq 0$, $u(x, \cdot)$ has at least a linear growth with respect to $y$ and hence $\int_0^\infty |u(x,y)|^p y^m\, dy=\infty$ (since $(m+1)/p>1-\alpha_2 >-1$). Then $g=0$ a.e. and  Proposition  \ref{Hardy in core}(ii) again gives $\|y^{\alpha_2 -1}D_yu\|_{L^p_m} \leq C\|y^{\alpha_2}D_{yy}u\|_{L^p_m}$.

If $\frac{m+1}{p} <1-\alpha_2$  we apply Proposition  \ref{Hardy in core} (i)  to $D_y u$ to deduce that  $\lim_{y \to 0}D_yu(x,y)=h(x)$  exists. If $h=0$, then Hardy inequality yields $y^{\alpha_2-1}D_y u \in L^p_m$. On the other hand,  $y^{\alpha_2-1}D_y u \in L^p_m$ implies $h=0$, since $y^{p(\alpha_2-1)}$ is not integrable with respect to the weight $y^m$.
\qed

\subsection{An alternative description of $W^{2,p}_{\mathcal N}(\alpha_1,\alpha_2,m)$}
We show an alternative description of $W^{2,p}_{\mathcal N}(\alpha_1,\alpha_2,m)$, adapted to the operator $D_{yy}+cy^{-1}D_y$.

	\begin{lem}\label{Lem Trace Dy in W}
		Let   $c\in\R$ and let us suppose that $\frac{m+1}{p}<c+1-\alpha_2$.  If $u \in  W^{2,p}_{loc}(\R^{N+1}_+)$ and $ u,\ y^{\alpha_2}\left(D_{yy}u+c\frac{D_yu}y\right) \in L^p_m$, then  the following properties hold.
		\begin{itemize}
			\item[(i)] The function  $v=y^cD_y u$ satisfies $v,D_yv\in L^1_{loc} \left(\R^{N}\times[0,\infty)\right)$  and therefore  has a   trace $v_0(x):=\lim_{y\to 0}y^c D_yu(x,y)\in L^1_{loc}(\R^N)$ at $y=0$. 
			\item[(ii)] $v_0=0$ if and only $y^{\alpha_2-1}D_y u \in L^p_m(\R^N \times [0,1])$. In this case
			\begin{align*}
				\left\|y^{\alpha_2-1}D_yu\right \|_{L^p_m}\leq C 	\left\|y^{\alpha_2}\left (D_{yy}u+cy^{-1}D_yu\right)\right\|_{L^p_m}
			\end{align*}
			with  $C=\left(c+1-\alpha_2-\frac{m+1}{p}\right)^{-1}>0$.
			\item[(iii)] If the stronger assumption $0<\frac{m+1}p\leq c-1$ holds then $v_0=0$ and $y^{\alpha_2-1}D_y u \in L^p_m(\R^N \times [0,1])$.
		\end{itemize}

	\end{lem}
	{\sc{Proof.}} Let   $v:=y^{c}D_yu$ and
	$$f:=y^{\alpha_2}\left(D_{yy}u+c\frac{D_yu}{y}\right)=y^{\alpha_2-c}D_yv\in L^p_m.$$ 
	Claim (i) is then a consequence of Proposition \ref{Hardy in core} (i) with $\beta=\alpha_2-c$.

	To prove (ii) we set $v_0(x)=\left(y^cD_yu\right)(x,0)$. Then one has  $g:=y^{\alpha_2-c-1}(v-v_0)\in L^p_m$ by Proposition \ref{Hardy in core} (ii) again. Then 
	$$y^{\alpha_2-1}D_yu=g+y^{\alpha_2-1-c}v_0$$ is  $L^p_m$-integrable near $y=0$ if and only if $v_0=0$, since  $\frac{m+1}{p} <c+1-\alpha_2$.
	
	Finally, when $v_0=0$, $y^{\alpha_2-1}=g=y^{\alpha_2-c-1}v$ and we can use  Proposition \ref{Hardy in core} (ii).
	
Let us  prove (iii).  Note that $c-1 <c+1-\alpha_2$, since $\alpha_2<2$. At the points where $v_0(x) \neq 0$, we have for  $0<y \leq \delta(x)$, $|D_yu(x, y)|\geq \frac 12|v_0(x)| y^{-c}$ which implies $|u(x,y)|\geq \frac 14|v_0(x)| y^{-c+1}$ for $0<y \leq \delta'(x)$, since $c>1$. This yields  $\int_0^\infty |u(x,y)|^p y^m\, dy=\infty$,  since $(m+1)/p\leq c-1$, and then $v_0=0$ a.e.
	\qed

	\medskip
	
	To provide an equivalent description of  $W^{2,p}_{\mathcal N}(\alpha_1, \alpha_2, m)$ we need the following simple lemma.
\begin{lem} \label{elliptic}
Assume that $u \in L^p(\R^N) \cap W^{2,1}_{loc}(\R^N)$ for some $1<p<\infty$ and that $\Delta u \in L^p(\R^N)$. Then $u \in W^{2,p}(\R^N)$.
\end{lem}
{\sc{Proof.}} Let $v \in W^{2,p}(\R^N)$ be such that $v-\Delta v=u-\Delta u$ and consider $w=u-v \in L^p(\R^N) \cap W^{2,1}_{loc}(\R^N)$. If $\phi \in C_c^\infty (\R^N)$, then
$$
0=\int_{\R^N}(w-\Delta w)\phi=\int_{\R^N}w(\phi-\Delta \phi).
$$
Since $w \in L^p(\R^N)$ the above identity extends by density to all $\phi \in W^{2,p'}(\R^N)$ and then, since $I-\Delta$ is invertible from $W^{2,p'}(\R^N)$ to $L^{p'}(\R^N)$, we have $\int_{\R^N} w g=0$ for every $g \in L^{p'}(\R^N)$, so that $w=0$ and $u=v \in W^{2,p}(\R^N)$.
\qed

We can now show an equivalent description of $W^{2,p}_{\mathcal N}(\alpha_1, \alpha_2, m)$, adapted to the degenerate operator $D_{yy}+cy^{-1}D_y$.

	\begin{prop}\label{Trace D_yu in W}
		Let   $c\in\R$ and $\frac{m+1}{p}<c+1-\alpha_2$.  Then
		\begin{align*} 
			W^{2,p}_{\mathcal N}(\alpha_1, \alpha_2, m)=&\left\{u \in  W^{2,p}_{loc}(\R^{N+1}_+): u,\ y^{\alpha_1}\Delta_xu\in L^p_m  \right. \\[1ex]			
			&\left.\hspace{10ex} y^{\alpha_2}\left(D_{yy}u+c\frac{D_yu}y\right) \in L^p_m\text{\;\;and\;\;}\lim_{y\to 0}y^c D_yu=0\right\}
		\end{align*}
and the norms $\|u\|_{W^{2,p}_{\mathcal N}(\alpha_1,\alpha_2,m)}$ and $$\|u\|_{L^p_m}+\|y^{\alpha_1}\Delta_x u\|_{L^p_m}+\|y^{\alpha_2}(D_{yy}u+cy^{-1}D_yu)\|_{L^p_m}$$ are equivalent on $W^{2,p}_{\mathcal N}(\alpha_1, \alpha_2, m)$.

		Finally,  when $0<\frac{m+1}p\leq c-1$ then 
				\begin{align*} 
			W^{2,p}_{\mathcal N}(\alpha_1, \alpha_2, m)=&\left\{u \in  W^{2,p}_{loc}(\R^{N+1}_+): u,\ y^{\alpha_1}\Delta_xu,   y^{\alpha_2}\left(D_{yy}u+c\frac{D_yu}y\right) \in L^p_m\right\}.
		\end{align*}
	\end{prop}
	{\sc{Proof.}} Let $\mathcal G$ be the space on the right hand side with the canonical norm indicated above. By Lemma \ref{Lem Trace Dy in W} $W^{2,p}_{\mathcal N}(\alpha_1, \alpha_2, m) \subset \mathcal G$ and the embedding is clearly continuous.

Conversely, let $u \in \mathcal G$. The estimate for $y^{\alpha_2-1}D_yu$ follows from  Lemma \ref{Lem Trace Dy in W}(ii)  and yields,  by difference, also that for $y^{\alpha_2}D_{yy}u$.   Since for $y\leq 1$ one has  $y^{\frac{\alpha_2}2}\leq y^{\alpha_2-1}$ it follows  that $y^{\frac{\alpha_2}2}D_yu\in L^p_m(\R^N \times [0,1])$ and  $y^{\frac{\alpha_2}2}D_yu\in L^p_m(\R^N \times [1,\infty])$ by Lemma \ref{inter1}.
	
	Finally, we prove the inequality
	$$\|y^{\frac{\alpha_1}2}D_xu\|_{L^p_m}+\|y^{\alpha_1}D_{x_ix_j}u\|_{L^p_m}\leq C\left (\|u\|_{L^p_m}+\|y^{\alpha_1}\Delta_x u\|_{L^p_m}\right ). $$
Since $u(\cdot,y) \in L^p(\R^N) \cap W_{loc}^{2,p}(\R^N)$ for a.e. $y>0$, the lemma above and the Calderon-Zygmund inequality yield
\begin{align*}
\int_{\R^N} |D_{x_i x_j}(x,y)|^p\,dx\leq C\int_{\R^N} |\Delta_x(x,y)|^p\,dx.
\end{align*}
Multiplying by $y^{p\alpha_1+m
}$ and integrating over $\R_+$ 
we obtain $\sum_{i,j=1}^n\|y^{\alpha_1} D_{x_ix_j}u\|_{L^p_m}\leq C\|y^{\alpha_1} \Delta_x u\|_{L^p_m}$. The estimate 
$$\|y^{\frac{\alpha_1}2} \nabla_{x}u\|_{L^p_m}\leq C\left(\|y^{\alpha_1} \Delta_x u\|_{L^p_m}+\|u\|_{L^p_m}\right)$$ can be obtained similarly using the  interpolative inequality
$$\|\nabla_x u(\cdot,y)\|_{L^p(\R^n)}\leq \epsilon \|\Delta_x u(\cdot,y)\|_{L^p(\R^n)}+\frac {C(N,p)} \epsilon \| u(\cdot,y)\|_{L^p(\R^n)}$$
 with  $\epsilon=y^{\frac{\alpha_1}2}$.
 
 The equality for $0<\frac{m+1}p\leq c-1$ follows from Lemma \ref{Lem Trace Dy in W}(iii).
\qed
		
		\medskip	
	We provide now another equivalent description of $W^{2,p}_{\mathcal N}(\alpha_1, \alpha_2, m)$ which involves a Dirichlet, rather than Neumann, boundary condition,  in a certain range of parameters.
		\begin{prop}\label{trace u in W op}
			Let   $c\geq 1$ and $\frac{m+1}{p}<c+1-\alpha_2$.  The following properties hold.
			\begin{itemize}
				\item[(i)] If $c>1$ then 
				\begin{align*}
					W^{2,p}_{\mathcal N}(\alpha_1, \alpha_2, m)=&\left\{u \in  W^{2,p}_{loc}(\R^{N+1}_+): u,\ y^{\alpha_1}\Delta_xu\in L^p_m,\right. \\[1ex]			
					&\left.\hspace{11ex}y^{\alpha_2}\left(D_{yy}u+c\frac{D_yu}y\right) \in L^p_m\text{\;and\;}\lim_{y\to 0}y^{c-1} u=0\right\}.
				\end{align*}
			\item[(ii)] If $c=1$ then 
\begin{align*}
	W^{2,p}_{\mathcal N}(\alpha_1, \alpha_2, m)=&\left\{u \in  W^{2,p}_{loc}(\R^{N+1}_+): u,\ y^{\alpha_1}\Delta_xu\in L^p_m,\right. \\[1ex]			
	&\left.\hspace{11ex}y^{\alpha_2}\left(D_{yy}u+c\frac{D_yu}y\right) \in L^p_m\text{\;and\;}\lim_{y\to 0} u(x,y)\in \C\right\}.
\end{align*}
			\end{itemize} 
		
		\end{prop}	
{\sc Proof. } Let us prove (i). By Proposition \ref{Trace D_yu in W} it is sufficient to show that the conditions $\lim_{y \to 0}y^cD_yu=0$ and $\lim_{y \to 0}y^{c-1}u=0$ are equivalent. We proceed as in Lemma \ref{Lem Trace Dy in W} setting   $v:=y^{c}D_yu$ and
	$$f:=y^{\alpha_2}\left(D_{yy}u+c\frac{D_yu}{y}\right)=y^{\alpha_2-c}D_yv\in L^p_m.$$ If $v_0(x)=\left(y^cD_yu\right)(x,0)$, then  $g:=y^{\alpha_2-c-1}(v-v_0)\in L^p_m$ by Proposition \ref{Hardy in core} (ii), and  
\begin{equation} \label{w1}	
D_yu=y^{1-\alpha_2}g+y^{-c}v_0.
\end{equation}
Then, since $c>1$, we can write for $0<y<1$
\begin{align}\label{eq1 trace u in W}
	u(x,1)- u(x,y)=\frac{1}{c-1}v_0(x)(y^{1-c}-1)+\int_y^1 s^{1-\alpha_2}g(x,s)\, ds
\end{align}
and
\begin{equation} \label{w2}
\int_y^1 s^{1-\alpha_2}|g(x,s)|\, ds \leq \|g\|_{L^p_m} \left(\int_y^1 s^{(1-\alpha_2 -\frac mp)p'} \right)^{\frac{1}{p'}} \leq C(1+y^{\gamma})
\end{equation}
with $\gamma=2-\alpha_2-(m+1)/p>1-c$ (when $\gamma=0$ the term $y^\gamma$ is substituted by $|\log y|^{\frac{1}{p'}}$). Since $c>1$, it follows that 
\begin{align*}
	\lim_{y \to 0} y^{c-1}u(x,y)=\frac{v_0(x)}{1-c}
\end{align*}
and therefore $\ds\lim_{y \to 0} y^{c-1}u(x,y)=0$ if and only if $v_0(x)=0$ or, by Lemma \ref{Lem Trace Dy in W}(ii), if $\ds\lim_{y \to 0}y^c D_yu(x,y)=0$.

To prove (ii) we proceed similarly. From \eqref{w1} with $c=1$ we obtain  
\begin{align*}
	u(x,1)- u(x,y)=-v_0(x)|\log y|+\int_y^1 s^{1-\alpha_2}g(x,s)\, ds,\qquad 0<y<1.
\end{align*}
The parameter $\gamma$ is positive, since $(m+1)/p<2-\alpha_2$ and the integral on the right hand side of \eqref{eq1 trace u in W} converges. 
Therefore $\ds\lim_{y \to 0} u(x,y)\in \C$ if and only if $v_0(x)=0$.
\qed

\begin{os}
	We point out that the function $v=y^{c-1}u$  above satisfies   $D_yv\in L^1 \left(Q\times[0,1]\right)$ for every cube $Q$. In particular 
 $D_y u\in L^1\left(Q\times [0,1]\right)$, if $\frac{m+1}{p} <2-\alpha_2$, by choosing $c=1$.

Indeed, if $c>1$,  using  \eqref{eq1 trace u in W}, \eqref{w2} with $v_0=0$ one has $y^{c-2}u\in L^1 \left(Q\times[0,1]\right)$. Then the equality $$D_{y}v=y^{c-1}D_yu+(c-1)y^{c-2}u=y^{c-\alpha_2}g+(c-1)y^{c-2}u$$
and $g \in L^p_m$ and H\"older inequality yield $y^{c-\alpha_2}g \in L^1(Q \times [0,1])$.

When $c=1$ then $v=u$ and we use \eqref{w1} with $v_0=0$ and then \eqref{w2}, since $\gamma>0$, as observed in the above proof.
\end{os}

\subsection{Approximation with smooth functions} \label{denso}

The main result of this section is a density property of smooth functions in $W^{2,p}_{\mathcal N}(\alpha_1,\alpha_2,m)$. We introduce the set
\begin{equation} \label{defC}
\mathcal{C}:=\left \{u \in C_c^\infty \left(\R^N\times[0, \infty)\right), \ D_y u(x,y)=0\  {\rm for} \ y \leq \delta\ {\rm  and \ some\ } \delta>0\right \}
\end{equation}
and its one dimensional version 
\begin{equation} \label {defD}
\mathcal{D}=\left \{u \in C_c^\infty ([0, \infty)), \ D_y u(y)=0\  {\rm for} \ y \leq \delta\ {\rm  and \ some\ } \delta>0\right \}.
\end{equation}
Let $$C_c^\infty (\R^{N})\otimes\mathcal D=\left\{u(x,y)=\sum_i u_i(x)v_i(y), \  u_i \in C_c^\infty (\R^N), \  v_i \in \cal D \right \}$$ (finite sums). Clearly $C_c^\infty (\R^{N})\otimes\mathcal D \subset \mathcal C$.
\begin{teo} \label{core gen}
If  $\frac{m+1}{p}>\alpha_1^-$
then $C_c^\infty (\R^{N})\otimes\mathcal D$ is dense in $W^{2,p}_{\mathcal N}(\alpha_1,\alpha_2,m)$.

\end{teo}

Note that the condition $(m+1)/p>\alpha_1^-$ or $m+1>0$ and $(m+1)/p+\alpha_1>0$ is necessary for the inclusion  $C_c^\infty (\R^{N})\otimes\mathcal D \subset W^{2,p}_{\mathcal N}(\alpha_1,\alpha_2,m)$. 
\medskip

For technical reason we start from the case $\alpha_2=0$ and write $\alpha$ for $\alpha_1$.
Then
$$W^{2,p}_{\mathcal N}(\alpha,0,m)=\left\{u\in W^{2,p}_{loc}(\R^{N+1}_+):\ u,\,  y^\alpha D_{x_ix_j}u,\ y^\frac{\alpha}{2} D_{x_i}u,\ D_{y}u\ D_{yy}u,\ \frac{D_yu}{y}\in L^p_m\right\}.$$

\medskip

We need some preliminary results which show the density of smooth functions with compact support in 
$W^{2,p}_{\mathcal N}(\alpha,0,m)$. In the first no restriction on $\alpha$ is needed.
\begin{lem} \label{supp-comp}
The functions in  $ W^{2,p}_{\mathcal N}(\alpha,0,m)$ having support in $\R^N \times [0,b[$ for some $b>0$ are dense in  $ W^{2,p}_{\mathcal N}(\alpha,0,m)$.
\end{lem}
{\sc Proof.}
Let $0\leq\phi\leq 1$ be a smooth function depending only on the $y$ variable which is equal to $1$ in $(0,1)$ and to $0$ for $y \ge 2$. 
Set $\phi_n(y)=\phi \left(\frac{y}{n}\right)$ and $u_n(x,y)=\phi_n(y)u(x,y)$.
Then $u_n\in W^{2,p}_{\mathcal N}(\alpha,0,m)$ and has compact support  in $\R^N \times [0,2n]$.  By dominated convergence $u_n \to u$ in $L^p_m$. Since $D_{x_ix_j}u_n=\phi_n  D_{x_ix_j}u$, $ D_{x_i}u_n=\phi_nD_{x_i}u$ we have $y^\alpha D_{x_ix_j}u_n\to y^\alpha D_{x_ix_j}u$, $y^\frac{\alpha}{2} D_{x_i}u_n\to y^\frac{\alpha}{2} D_{x_i}u$ , by dominated convergence again.

For the convergence of the $y$-derivatives, we observe that 
$|D_y \phi_n|\leq \frac{C}{n}\chi_{[n,2n]}$, $|D_{yy} \phi_n| \leq \frac{C}{n^2}\chi_{[n,2n]}$. Since  
$D_y u_n =\phi_n D_y u+ D_y \phi_n u$ and $D_{yy} u_n =\phi_n D_{yy} u+2D_y\phi_n D_y u+ uD_{yy}\phi_n$, we have  also $D_y u_n \to D_y u$, 
$D_{yy}u_n\to  D_{yy}u$ and   $\frac{D_yu_n}{y}\to \frac{D_yu}{y}$ in $L^p_m$. 
\qed

\begin{lem} \label{supp-comp-x}
Assume that $\frac{m+1}{p}<2$ and $\frac{m+1}{p}+\alpha>0$. Then the functions in  $ W^{2,p}_{\mathcal N}(\alpha,0,m)$ with compact support  are dense in  $ W^{2,p}_{\mathcal N}(\alpha,0,m)$.
\end{lem}
{\sc Proof.} Let $u\in W^{2,p}_{\mathcal N}(\alpha,0,m)$.  By Lemma \ref{supp-comp}, we may assume that $u$ has support in  $\R^N \times [0,b[$ for some $b>0$. 
Let $0\leq\phi\leq 1$ be a smooth function depending only on the $x$ variable which is equal to $1$ if $|x|\leq 1$ and to $0$ for $|x| \ge 2$. 
Set $\phi_n(x)=\phi \left(\frac{x}{n}\right)$ and $u_n(x,y)=\phi_n(x)u(x,y)$.
Then $u_n\in W^{2,p}_{\mathcal N}(\alpha,0,m)$ and has compact support.  By dominated convergence $u_n \to u$ in $L^p_m$. Moreover, since $D_y u_n =\phi_n D_y u$ and $D_{yy} u_n =\phi_n D_{yy} u$, we have  immediately $D_{yy}u_n\to  D_{yy}u$,  $\frac{D_yu_n}{y}\to \frac{D_yu}{y}$ in $L^p_m$. 

Concerning the derivatives with respect to the $x$ variable, we have
$|D_{x_i} \phi_n(x)|\leq \frac{C}{n}\chi_{[n,2n]}(|x|)$, $|D_{x_ix_j } \phi_n(x)|\leq \frac{C}{n^2}\chi_{[n,2n]}(|x|)$ and
\begin{align}\label{eq1 lem supp}
\nonumber D_{x_i}u_n&=\phi_n D_{x_i}u+u D_{x_i}\phi_n, \\
D_{x_jx_i}u_n&=\phi_n D_{x_ix_j}u+D_{x_j}\phi_n D_{x_i}u+D_{x_j}u D_{x_i}\phi_n +u  D_{x_ix_j}\phi_n. 
\end{align}

Let us show that $y^\alpha u,\ y^{\frac{\alpha}2} u\in L^p_m$.  Since $u$ has  support in $\R^N \times [0,b[$ this is  trivial for $\alpha\geq 0$. When $\alpha<0$ let   $f(x,y)=\frac{u(x,y)-u(x,0)}{y^2}$ so that
\begin{align*}
y^\alpha u=y^{\alpha+2}f+y^\alpha  u(\cdot,0).
\end{align*}
 By  Proposition  \ref{Hardy in core}, $f\in L^p_m$ and $u(\cdot,0)\in L^p(\R^N)$.   Since $u$ and $f$ have support in $\R^N \times [0,b[$,  the  assumption  $-\alpha<\frac{m+1}{p}<2$ then implies that $y^\alpha u\in L^p_m$ and also  $y^{\frac{\alpha}2} u\in L^p_m$.
Using the classic interpolative inequality
$$\|\nabla_x u(\cdot,y)\|_{L^p(\R^N)}\leq \epsilon \|\Delta_x u(\cdot,y)\|_{L^p(\R^N)}+\frac {C(N,p)} \epsilon \| u(\cdot,y)\|_{L^p(\R^N)}$$
 with  $\epsilon=1$ we easily get (after raising to the power $p$,  multiplying by $y^\alpha$ and integrating in $y$), $y^\alpha \nabla_x u\in L^p_m$.
Using \eqref{eq1 lem supp} and the fact that $y^\alpha  u,\ y^{\frac{\alpha}2} u,\ y^\alpha \nabla_x u\in L^p_m$
we deduce using dominated convergence that
$y^\alpha D_{x_ix_j}u_n\to y^\alpha D_{x_ix_j}u$, $y^\frac{\alpha}{2} D_{x_i}u_n\to y^\frac{\alpha}{2} D_{x_i}u$ in $L^p_m$. 
\qed

In the next lemma we add regularity with respect to the $x$-variable.

\begin{lem}
\label{smooth}
Let  $\frac{m+1}{p}<2$, $\frac{m+1}{p}+\alpha>0$ and  $u\in W^{2,p}_{\mathcal N}(\alpha,0,m)$ with compact support. Then there exist $(u_n)_{n\in\N}\subset W^{2,p}_{\mathcal{N}}(\alpha,0,m)$ with compact  support  such that $u_n$ converges to $u$ in $W^{2,p}_{\mathcal{N}}(\alpha,0,m)$ and, for every $y\geq 0$,  $u_n(\cdot,y)$ belongs to $ C^\infty(\R^N)$ and has bounded $x$-derivatives of any order.
\end{lem}
{\sc Proof.} Let $u\in W^{2,p}_{\mathcal{N}}(\alpha,0,m)$ as above and let us fix a standard sequence of mollifiers in $\R^N$  $\rho_n=n^N\rho(nx)$ where $0\leq \rho\in C_c^\infty(\R^N)$, $\int_{\R^N}\rho(x)\ dx=1$. Let us  set $u_n(x,y)=\left(\rho_n\ast u(\cdot,y)\right)(x)$ where $\ast$ means convolution with respect to the $x$ variable.

By Lemma \ref{L1} and Proposition \ref{Hardy in core},  $u(\cdot, y)\in L^p(\R^N)$ for every $y\geq 0$ and
therefore, by standard properties,   $u_n$  has a compact support and $u_n(\cdot,y)\in C^\infty_b(\R^N)$ for every $y\geq 0$. By Young's inequality
\begin{align*}
\|u_n(\cdot,y)\|_{L^p(\R^N)}\leq \|u(\cdot,y)\|_{L^p(\R^N)},\qquad u_n(\cdot,y)\to u(\cdot,y)\quad\text{in}\quad L^p(\R^N),\qquad \forall y\geq 0.
\end{align*}
Raising to the power $p$, multiplying by $y^m$ and by integrating with respect to $y$, we get
\begin{align*}
\|u_n\|_{L^p_m}\leq \|u\|_{L^p_m}
\end{align*} 
 which, using dominated convergence with respect to $y$,   implies $u_n\to u$ in $L^p_m$.
 
Using the equalities 
\begin{align*}
y^\alpha D_{x_ix_j}u_n&=\rho_n\ast (y^\alpha D_{x_ix_j}u),\qquad y^\frac{\alpha}{2} D_{x_i}u_n=\rho_n\ast (y^\frac{\alpha}{2} D_{x_i}u),\\[1ex]
D_{yy}u_n&=\rho_n\ast D_{yy}u,\hspace{14.5ex} y^\gamma D_{y}u_n=\rho_n\ast (y^\gamma D_{y}u),
\end{align*}
$\gamma=0,1$ a similar argument as before yields  $u_n\to u$ in $W^{2,p}_n(\alpha,0,m)$.\\\qed

We can now prove a weaker version of  Theorem \ref{core gen} when $\alpha_2=0$.
\begin{prop} \label{corend}
If  $\frac{m+1}{p}>\alpha^-$ then  $\mathcal C$, defined in \eqref{defC},  is dense in $W^{2,p}_{\mathcal{N}}(\alpha,0,m)$.
\end{prop}
{\sc Proof.}  (i) We first consider the case $\frac{m+1}{p}>2$.  Let $u\in W^{2,p}_{\mathcal{N}}(\alpha,0,m)$ which, by  Lemma \ref{supp-comp}, we may assume to have the  support in $\R^N\times [0,b]$. Let $\phi$ be a smooth function depending only on $y$,  equal to $0$ in $(0,1)$ and to $1$ for $y \ge 2$. Let $\phi_n(y)=\phi (ny)$ and $u_n(x,y)=\phi_n(y)u(x,y)$.  Then
\begin{align*}
D_{x_ix_j}u_n &=\phi_n D_{x_ix_j}u,\hspace{12ex}D_{x_i}u_n =\phi_n D_{x_i}u,\\[1ex]
D_y u_n &=\phi_n D_y u+D_y\phi_n u,\qquad D_{yy} u_n =\phi_n D_{yy} u+2D_y\phi_n D_yu+ uD_{yy}\phi_n. 
\end{align*}
By dominated convergence $u_n \to u$, $y^\alpha D_{x_ix_j}u_n \to y^\alpha D_{x_ix_j}u$, $y^\frac{\alpha}{2} D_{x_i}u_n \to y^\frac{\alpha}{2}  D_{x_i}u$ in $L^p_m$.
Let us consider now the terms containing the $y$ derivatives  and observe that 
\begin{align}\label{sti cut 2}
|D_{y} \phi_n|\leq Cn\chi_{[\frac 1 n,\frac 2{n}]}\leq \frac{2C}{y}\chi_{[\frac 1 n,\frac 2{n}]},\qquad |D_{yy } \phi_n|\leq C n^2\chi_{[\frac 1 n,\frac 2{n}]}\leq \frac{4C}{y^2}\chi_{[\frac 1 n,\frac 2{n}]}.
\end{align}
Using these estimates and since $y^{-2}u\in L^p_m $ by Proposition  \ref{Hardy in core}

$$\frac{D_y u_n}{y} =\phi_n \frac{D_y u}{y}+\frac{u}{y} (D_y\phi_n) \to \frac{D_y u}{y}$$
in $L^p_m$, by dominated convergence.

In a similar way one shows that $D_y u_n \to D_yu$ and $D_{yy}u_n \to D_{yy}u$ in $L^p_m$ and hence functions with compact support in $\R^n\times ]0,\infty[$ are dense in $W^{2,p}_{\mathcal{N}}(\alpha,0,m)$. At this point, a standard smoothing by convolutions shows the density of $C_c^\infty (\R^N \times ]0,\infty[)$ in $W^{2,p}_{\mathcal{N}}(\alpha,0,m)$.


(ii) Let  $\frac{m+1}{p}=2$. We proceed similarly to (i) and fix $u\in W^{2,p}_{\mathcal{N}}(\alpha,0,m)$ with  support in $\R^N\times [0,b]$. 
Let $\phi$ be a smooth function which is equal to $0$ in $\left(0,\frac{1}{4}\right)$ and to $1$ for $y \ge \frac{1}{2}$. 
Let $\phi_n(y)=\phi\left(y^\frac{1}{n}\right)$ and $u_n=\phi_n u$. By dominated convergence it is immediate to see that  $u_n \to u$,\; $y^{\alpha}D_{x_ix_j}u_n\to y^{\alpha}D_{x_ix_j}u$,\; $y^\frac{\alpha}{2}D_{x_i}u_n\to y^\frac{\alpha}{2}D_{x_i}u$  in $L^p_m$. To treat the terms concerning the $y$ derivatives we observe that 
\begin{align}\label{beh cut}
\nonumber |\phi_n'|&=\left|\frac{1}{n}\phi'\left(y^\frac{1}{n}\right)y^{\frac{1}{n}-1} \right|\leq \frac{C}{ny}\chi_{[(\frac 1 4)^n,(\frac 1{2})^n]}\\[1ex]
 | \phi_n''|&=\left|\frac{1}{n^2}\phi''\left(y^\frac{1}{n}\right)y^{\frac{2}{n}-2}+\frac 1 n(\frac 1 n-1)\phi'\left(y^\frac{1}{n}\right)y^{\frac{1}{n}-2}\right|\leq  \frac{C}{ny^2}\chi_{[(\frac 1 4)^n,(\frac 1{2})^n]}.
\end{align}

Moreover, 
\begin{align*}
D_y u_n =\phi_n  D_y u+\phi_n'u,\qquad D_{yy} u_n=\phi_n D_{yy} u+2\phi_n'D_yu+\phi_n''u.
\end{align*}
Then $\frac 1 y D_y u_n\to \frac 1 y D_y u$ in $L^p_m$ since $\phi_n  \frac{D_y u}{y}\to \frac 1 y D_y u$  by dominated convergence and $\phi_n'\frac{u}{y}\to 0$. In fact,   using  \eqref{beh cut} and \eqref{behaviour}  of Lemma \ref{int-uMaggiore} we have
\begin{align*}
\left\|\phi_n'\frac{u}{y}\right\|^p_{L^p_m}\leq \frac{C}{n^p}\int_{(\frac 1 4)^n}^{(\frac 1 2)^n} |\log y|^{p-1}y^{m-2p}\,dy=\frac{C}{n^{2p}}\int_{\frac 1 4}^{\frac 1 2} |\log s|^{p-1}s^{-1}\,dy
\end{align*}
which tends to $0$ as $n\to\infty$. 

Concerning the second order derivative we  have $ D_{yy} u_n\to D_{yy} u$ since  $\phi_n D_{yy} u\to D_{yy} u$  by dominated convergence and  the other terms tend to $0$. Indeed proceeding as before we have $|\phi_n'D_yu|\leq C\frac{C}{n}\chi_{[(\frac 1 4)^n,(\frac 1{2})^n]}\frac{|D_y u|}{y}$ which tends to $0$ by dominated convergence. Finally, 
 \begin{align*}
\|\phi_n''u\|^p_{L^p_m}\leq \frac{C}{n^p}\int_{(\frac 1 4)^n}^{(\frac 1 2)^n} |\log y|^{p-1}y^{m-2p}\,dy=\frac{C}{n^{2p}}\int_{\frac 1 4}^{\frac 1 2} |\log s|^{p-1}s^{-1}\,dy
\end{align*}
which tends to $0$ as $n\to\infty$.

Now the proof is as for (i) and shows that $C_c^\infty (\R^N \times ]0,\infty[)$ is dense in $W^{2,p}_{\mathcal{N}}(\alpha,0,m)$.

(iii) Let assume finally that $\frac{m+1}{p}<2$. By Lemmas \ref{supp-comp-x}, \ref{smooth} we may assume that $u$ has compact support and that for every $y \geq 0$, $u(\cdot,y) \in C^\infty_b(\R^N)$.

By Proposition \ref{Hardy in core},   $\frac{u-u(\cdot,0)}{y^2}\in L^p_m$. Let $\phi$ be a smooth function  equal to $0$ in $(0,1)$ and to $1$ for $y \ge 2$ and $\phi_n(y)=\phi (ny)$. Setting 
$$u_n(x,y)=(1-\phi_n(y))u(x,0)+\phi_n(y)u(x,y),$$
then  
\begin{align*}
D_{x_i}u_n &=(1-\phi_n)D_{x_i} u(\cdot,0)+ \phi_n D_{x_i}u,\\[1ex]
D_{x_ix_j}u_n &=(1-\phi_n)D_{x_ix_j} u(\cdot,0)+ \phi_nD_{x_ix_j}u,\\[1ex]
D_y u_n &=\phi_n' (u-u(\cdot,0))+\phi_nD_{y}u,\\[1ex]
D_{yy} u_n &=\phi_n'' (u-u(\cdot,0))+2\phi_n'D_yu+\phi_nD_{yy}u.
\end{align*}

It follows that $u_n \to u$,\; $y^\alpha D_{x_ix_j}u_n \to y^\alpha D_{x_ix_j}u$, \; $y^\frac{\alpha}{2} D_{x_i}u_n \to y^\frac{\alpha}{2}  D_{x_i }u$ in $L^p_m$. Since the argument is always the same, let us explain it for $u_n$. The term $\phi_n u$ converges to $u$ by dominated convergence and $(1-\phi_n)u(\cdot, 0)$ converges to zero since $u(\cdot,0)$ is bounded with compact support.

Using \eqref{sti cut 2} one has 
 $$\frac{|\phi_n' (u-u(\cdot,0))|}{y}\leq C \chi_{[\frac{1}{n},\frac{2}{n}]}(y)\frac{|u-u(\cdot,0)|}{y^2}$$
which tend to $0$  in $L^p_m$ by dominated convergence and then  $\frac{D_y u_n}{y}$ converges to $\frac{D_{y}u}{y}$ in $L^p_m$. 
Similarly $D_{yy} u_n$ converges $D_{yy}u$ in $L^p_m$.

Each function $u_n$ has compact support, does not depend on $y$ for small $y$ and is smooth with respect to the $x$ variable for any fixed $y$. Smoothness with respect to $y$ is however not yet guaranteed. This last property can be added by taking appropriate convolutions in $y$ with a compact support mollifier.
\qed

We can now prove the  general density result.

\medskip

{\sc{Proof of Theorem \ref{core gen}} }
The density of $\mathcal C$, defined in \eqref{defC},  in  $W^{2,p}_{\mathcal N}(\alpha_1,\alpha_2,m)$  follows by  Lemma \ref{Sobolev eq} and  Proposition \ref{corend} since  the isometry $T_{0,-\frac{\alpha_2}2}$ isometrically maps  dense subsets of  $W^{2,p}_{\mathcal{N}}(\tilde \alpha,0,\tilde m)$ into dense subsets of $W^{2,p}_{\mathcal{N}}(\alpha_1,\alpha_2,m)$ and, since  $\alpha_2<2$, leaves invariant $\mathcal{C}$. Note also that the conditions $(m+1)/p>\alpha_1^-$ and $(\tilde m+1)/p>\tilde\alpha^-$ are equivalent, since $\alpha_2<2$, again.

In order to prove the density of $C_c^\infty (\R^{N})\otimes\mathcal D$, we may therefore assume that $u$ is in $\mathcal C$, that is $u \in C_c^\infty (\R^{N} \times [0, \infty))$ and $ D_y u(x,y)=0$ for $y \leq \delta$ for some $\delta>0$.
Let  $\eta$ be a smooth function depending only on the $y$ variable which is equal to $1$ in $[0,\frac{\delta}{2}]$ and to $0$ for $y \ge \delta$. Then, since $D_y u(x,y)=0$ for  $y \leq \delta$,
$$u(x,y)=\eta (y)u(x,y)+(1-\eta(y))u(x,y)=\eta (y)w(x)+(1-\eta(y))u(x,y)=u_1(x,y)+u_2(x,y)$$
with $u_1(x,y)=\eta(y)w(x)$, $w(x)$ depending only on the $x$ variable.
Observe now that  $u_2(x,y)=(1-\eta(y))u(x,y)=0$ in $[0,\frac{\delta}{2}]$ and outside the support of $u$, therefore it belongs to  $C^\infty_c(\R^{N+1}_+)$ and the approximation
with respect to the  $W^{2,p}(\R^{N+1}_+)$ norm by functions in   $C_c^\infty (\R^{N})\otimes C_c^\infty (]0, \infty[)$ is standard (just use a sequence of polynomials converging uniformly to $u_2$ with all first and second order derivatives on a cube containing the support of $u_2$ and truncate outside the cube by a cut-off of the same type).  This proves the result.\qed

\begin{os}
From 
the proofs of Proposition \ref{corend} and Theorem \ref{core gen} it follows that  if $u\in W^{2,p}_{\mathcal N}(\alpha_1,\alpha_2,m)$ has  support in $\R^N\times[0,b]$, then there exists a sequence $\left(u_n\right)_{n\in\N}\in\mathcal C$  such that $ \mbox{supp }u_n\subseteq \R^N\times[0,b]$  and  $u_n\to u$ in $W^{2,p}_{\mathcal N}(\alpha_1,\alpha_2,m)$.
\end{os}

\begin{cor}
	\label{Core C c infty}
	Assume $\frac{m+1}{p}\geq 2-\alpha_2$ and $\frac{m+1}{p}>\alpha_1^-$. Then  $ C_c^\infty (\R^{N+1}_+) $  and 
	$C_c^\infty (\R^{N})\otimes C_c^\infty \left(]0, \infty[\right)$ are dense in $W^{2,p}_{\mathcal N}(\alpha_1,\alpha_2,m)$. 
\end{cor}
{\sc Proof. } This follows from the 
 the proofs of Proposition \ref{corend} and of Theorem \ref{core gen}. 
\qed
	
\medskip

 Specializing  Proposition \ref{Hardy in core} to $W^{2,p}_{\mathcal N}(\alpha_1, \alpha_2, m)$ we  get  the following corollary.

\begin{cor}\label{Hardy Rellich Sob}
	Let   $\frac{m+1}{p}>\alpha_1^-$. The following properties hold for any $u\in W^{2,p}_{\mathcal N}(\alpha_1, \alpha_2, m)$.
	\begin{itemize}
		\item[(i)] If $\frac{m+1}p>1-\frac{\alpha_2}2$ then
		\begin{align*}
			\|y^{\frac{\alpha_2}2-1}u\|_{L^p_m}\leq C \|y^{\frac{\alpha_2}2}D_{y}u\|_{L^p_m}.
		\end{align*}
		\item[(ii)] If $\frac{m+1}p>2-\alpha_2$ then
		\begin{align*}
			\|y^{\alpha_2-2}u\|_{L^p_m}\leq C \|y^{\alpha_2-1}D_{y}u\|_{L^p_m}.
		\end{align*}
		\item[(iii)] If $\frac{m+1}p<2-\alpha_2$ then
		\begin{align*}
			\|y^{\alpha_2-2}(u-u(\cdot,0))\|_{L^p_m}\leq C \|y^{\alpha_2-1}D_{y}u\|_{L^p_m}.
		\end{align*}
		
		\item[(iv)] If  $\alpha_2-\alpha_1<2$ and  $\frac{m+1}p>1-\frac{\alpha_1+\alpha_2}{2}$, $\frac{m+1}p>\alpha_1^-$ then
		\begin{align*}
			\|y^{\frac{\alpha_1+\alpha_2}{2}-1}\nabla_{x}u\|_{L^p_m}\leq C \|y^\frac{\alpha_1+\alpha_2}{2} D_{y}\nabla_x u\|_{L^p_m}.
		\end{align*}
	\end{itemize}
\end{cor}
{\sc{Proof.}}  By density we may assume that $u \in C_c^\infty (\R^{N})\otimes\mathcal D$. All points follow by applying Proposition \ref{Hardy in core}  to  $u$ in the cases (i), (ii) and (iii) and to $\nabla_x u$ in the case  (iv), recalling  Proposition \ref{Sec sob derivata mista}.\qed
\medskip

\section{The space $W^{2,p}_{\mathcal R}(\alpha_1, \alpha_2, m)$} \label{Sec sob  min domain}

We  consider an integral version of Dirichlet boundary conditions, namely  a weighted summability of $y^{-2}u$ and introduce for $m \in \R$, $\alpha_2<2$
\begin{equation} \label{dominiodirichlet}
	W^{2,p}_{\mathcal R}(\alpha_1, \alpha_2, m)=\{u \in  W^{2,p}(\alpha_1, \alpha_2, m): y^{\alpha_2-2}u \in L^p_m\}
\end{equation}
with the norm $$\|u\|_{W^{2,p}_{\mathcal R}(\alpha_1, \alpha_2, m)}=\|u\|_{W^{2,p}(\alpha_1, \alpha_2, m)}+\|y^{\alpha_2-2}u\|_{L^p_m}.$$
We remark that $W^{2,p}_{\mathcal R}(\alpha_1, \alpha_2, m)$ will be considered for every $m \in \R$ whereas $W^{2,p}_{\mathcal N}(\alpha_1, \alpha_2, m)$ only for $m+1>0$. The symbol $\mathcal R$ stands for "Rellich", since Rellich inequalities concern with the summability of $y^{-2}u$.

\begin{prop} \label{RN} The following properties hold.
	\begin{itemize}
		\item[(i)]  if $u \in W^{2,p}_{\mathcal R}(\alpha_1, \alpha_2, m)$ then $y^{\alpha_2-1}D_y u \in L^p_m.$
		\item[(ii)] If $\alpha_2-\alpha_1<2$ and $\frac{m+1}{p}>2-\alpha_2$, then $W^{2,p}_{\mathcal R}(\alpha_1, \alpha_2, m) =  W^{2,p}_{\mathcal N}(\alpha_1, \alpha_2, m)=W^{2,p}(\alpha_1, \alpha_2, m)$, with equivalence of the corresponding norms. In particular, $C_c^\infty (\R^{N+1}_+)$ is dense in $W^{2,p}_{\mathcal R}(\alpha_1, \alpha_2, m) $.
		

	\end{itemize}
\end{prop}
{\sc Proof. } The proof of (i) follows by integrating with respect to $x$ the inequality of Lemma \ref{inter}. The proof of  (ii) follows from  Proposition \ref{neumann}(i) and Corollary \ref{Hardy Rellich Sob}(ii), after noticing that $\alpha_2-\alpha_1 <2$ and $\frac{m+1}{p}>2-\alpha_2$ yield  $\frac{m+1}{p}>\alpha_1^-$. The density of $C_c^\infty (\R^{N+1}_+)$ in $W^{2,p}_{\mathcal R}(\alpha_1, \alpha_2, m)$ now follows from Corollary \ref{Core C c infty}.
\qed

Finally, we investigate the action of the multiplication operator $T_{k,0}:u\mapsto y^ku$. The following lemma  is the companion of Lemma \ref{Sobolev eq} which deals with the transformation $T_{0,\beta}$.
\begin{lem}\label{isometryRN}
	\label{y^k W}
	Let   $\alpha_2-\alpha_1<2$ and  $\frac{m+1}{p}>2-\alpha_2$. For every $k\in\R$
	\begin{align*}
		T_{k,0}:  W^{2,p}_{\mathcal N}(\alpha_1, \alpha_2, m) \to  W^{2,p}_{\mathcal R}(\alpha_1, \alpha_2, m-kp)
	\end{align*}
	
	is an isomorphism (we shall write $y^k  W^{2,p}_{\mathcal N}(\alpha_1, \alpha_2, m)= W^{2,p}_{\mathcal R}(\alpha_1, \alpha_2, m-kp)$).
\end{lem}
{\sc{Proof.}} Let   $u=y^{k}v$ with $v\in  W^{2,p}_{\mathcal N}(\alpha_1, \alpha_2, m)$. Since all $x$-derivatives commute with $T_{k,0}$ we deal only with the $y$-derivatives. We observe that 
\begin{align*}
	D_yu=y^k(D_yv+k\frac{v}y),\qquad
	D_{yy}u=y^k\left(D_{yy}v+2k \frac{D_y v}{y}+ k(k-1)\frac v{y^2}\right).
\end{align*}
Corollary \ref{Hardy Rellich Sob} yields 
\begin{align*}
	\|y^{\alpha_2-2}v\|_{L^p_m}+\|y^{\frac{\alpha_2}2-1}v\|_{L^p_m}+\|y^{\alpha_2-1}D_y v\|_{L^p_m}\leq c \|v\|_{W^{2,p}_{\mathcal N}(\alpha_1, \alpha_2, m)}
\end{align*}
and then $u \in W^{2,p}_{\mathcal R}(\alpha_1, \alpha_2, m-kp)$. Conversely, if $u \in W^{2,p}_{\mathcal R}(\alpha_1, \alpha_2, m-kp)$, then $y^{\alpha_2-1}D_y u \in L^p_{m-kp}$ by Proposition \ref{RN}(i) and similar formulas as above show that $y^{\alpha_2-1}D_y v, y^{\alpha_2} D_{yy}v \in L^p_m$. Since $y^{\alpha_2/2-1} \leq 1+y^{\alpha_2-2}$, then $y^{\alpha_2/2-1} u \in L^p_{m-kp}$ and  $y^{\alpha/2}D_y v \in L^p_m$.
\qed

\bibliography{../TexBibliografiaUnica/References}

\begin{thebibliography}{10}

\bibitem{met-calv-negro-spina}
{\sc {Calvaruso}, G., {Metafune}, G., {Negro}, L., and {Spina}, C.}
\newblock Optimal kernel estimates for elliptic operators with second order
  discontinuous coefficients.
\newblock {\em Journal of Mathematical Analysis and Applications 485}, 1
  (2020), 123763.

\bibitem{Geymonat-Grisvard}
{\sc {G. Geymonat, P. Grisvard}}.
\newblock In {\em Problemi ai limiti ellittici negli spazi di Sobolev con
  peso}. Cooperativa libraria universitaria pavese, 1965, pp.~1--57.

\bibitem{grisvard}
{\sc Grisvard, P.}
\newblock Espaces interm{\'e}diaires entre espaces de sobolev avec poids.
\newblock {\em Annali della Scuola Normale Superiore di Pisa - Classe di
  Scienze 3e s{\'e}rie, 17}, 3 (1963), 255--296.

\bibitem{met-negro-soba-spina}
{\sc {Metafune}, G., {Negro}, L., {Sobajima}, M., and {Spina}, C.}
\newblock {Rellich} inequalities in bounded domains.
\newblock {\em Mathematische Annalen 379\/} (2021), 765--824.

\bibitem{MNS-Sharp}
{\sc {Metafune}, G., {Negro}, L., and {Spina}, C.}
\newblock Sharp kernel estimates for elliptic operators with second-order
  discontinuous coefficients.
\newblock {\em Journal of Evolution Equations 18\/} (2018), 467--514.

\bibitem{MNS-Grad}
{\sc {Metafune}, G., {Negro}, L., and {Spina}, C.}
\newblock Gradient estimates for elliptic operators with second-order
  discontinuous coefficients.
\newblock {\em Mediterranean Journal of Mathematics 16}, 138 (2019).

\bibitem{MNS-Grushin}
{\sc {Metafune}, G., {Negro}, L., and {Spina}, C.}
\newblock {$L^p$} estimates for {Baouendi–Grushin} operators.
\newblock {\em Pure and Applied Analysis 2}, 3 (2020), 603--625.

\bibitem{MNS-Max-Reg}
{\sc {Metafune}, G., {Negro}, L., and {Spina}, C.}
\newblock Maximal regularity for elliptic operators with second-order
  discontinuous coefficients.
\newblock {\em Journal of Evolution Equations\/} (2020).

\bibitem{MNS-Caffarelli}
{\sc {Metafune}, G., {Negro}, L., and {Spina}, C.}
\newblock {$L^p$} estimates for the {Caffarelli-Silvestre} extension operators.
\newblock {\em Submitted\/} (2021).
\newblock Online preprint: \url{https://arxiv.org/abs/2103.10314v1}.

\bibitem{MNS-PerturbedBessel}
{\sc {Metafune}, G., {Negro}, L., and {Spina}, C.}
\newblock Degenerate operators on the half-line.
\newblock {\em Submitted\/} (2022).
\newblock Online preprint on arxiv.

\bibitem{MNS-CompleteDegenerate}
{\sc {Metafune}, G., {Negro}, L., and {Spina}, C.}
\newblock A unified approach to degenerate problems in the half-space.
\newblock {\em Submitted\/} (2022).
\newblock Online preprint on arxiv.

\bibitem{met-soba-spi-Rellich}
{\sc Metafune, G., Sobajima, M., and Spina, C.}
\newblock Weighted {Calder{\'o}n}--{Zygmund} and {Rellich} inequalities in
  {$L^p$}.
\newblock {\em Mathematische Annalen 361}, 1 (Feb 2015), 313--366.

\bibitem{morel}
{\sc Morel, H.}
\newblock Introduction de poids dans l'\'etude de probl\`emes aux limites.
\newblock {\em Annales de l'Institut {Fourier} 12\/} (1962), 299--413.

\bibitem{necas}
{\sc {Necas}, J.}
\newblock {\em Direct methods in the theory of elliptic equations}.
\newblock Springer, Berlin, 2012.

\bibitem{Negro-Spina-Asympt}
{\sc Negro, L., and Spina, C.}
\newblock Asymptotic behaviour for elliptic operators with second-order
  discontinuous coefficients.
\newblock {\em Forum Mathematicum 32}, 2 (2020), 399--415.

\end{thebibliography}

\end{document}